\documentclass[10pt,twocolumn,twoside]{IEEEtran}
\usepackage{amssymb}
\usepackage{amsmath}\usepackage{graphicx}
\usepackage{times}\hfuzz=10pt 

\usepackage{color}

\newcommand{\rd}[1]{\color{red}#1\color{black}} 
 
\newtheorem{theorem}{Theorem}
\newtheorem{lemma}{Lemma}
\newtheorem{definition}{Definition}

\def\e{\varepsilon}

\def\defi{\stackrel{{\scriptscriptstyle \Delta}}{=}}

\def\a{\alpha}

\def\o{\omega}
\def\O{\Omega}

\def\Y{{\cal Y}}

\def\w{\widehat}
\def\Ind{{\mathbb{I}}}

\def\esssup{\mathop{\rm ess\, sup}}

\def\Re{{\rm Re\,}}

\def\R{{\bf R}}

\def\Z{{\cal Z}}
\def\ZZ{{\bf Z}}

\def\g{\gamma}
\def\C{{\bf C}}

\def\X{{\cal X}}

\def\oo{\bar}

\def\U{{\cal U}}

\def\T{{\mathbb{T}}}

\def\k{\kappa}

\newcommand{\be}{\begin{equation}}
\newcommand{\ee}{\end{equation}}
\newcommand{\bd}{\begin{displaymath}}
\newcommand{\ed}{\end{displaymath}}
\newcommand{\ba}{\begin{array}{ll}}
\newcommand{\ea}{\end{array}}
\newcommand{\baa}{\begin{eqnarray}}
\newcommand{\eaa}{\end{eqnarray}}
\newcommand{\baaa}{\begin{eqnarray*}}
\newcommand{\eaaa}{\end{eqnarray*}}

\def\oo{\bar}

\def\a{\alpha}

\def\K{{\cal K}}
\def\Ko{K}
\def\ko{k}
\def\Yo{Y}
\def\yo{y}
\def\ew{\left(e^{i\o}\right)}
\date{Submitted: December 7, 2011. Revised: August 21, 2018 }
\title{
Predictors for discrete time processes with energy decay on higher
frequencies }
\author{
Nikolai Dokuchaev}
 
\begin{document}
\def\break{}%
\def\brea{}
\def\breakk{}
\def\brea{\nonumber\\ }\def\breakk{\nonumber\\&&} 
 \vspace{-0.5cm}      \maketitle
\begin{abstract} The predictability of
discrete-time processes is studied in a deterministic setting. A
family of one-step-ahead predictors is suggested for processes of
which the energy decays at higher frequencies. For such processes,
the prediction error can be made arbitrarily small. The predictions
can be robust with respect to the noise contamination at higher
frequencies. \let\thefootnote\relax\footnote{This is a
pre-copy-editing, author-produced PDF of an article published in 
{\em IEEE Transactions on Signal Processing}, November 2012   {\bf 60}, No. 11, 6027-6030.The 
publisher-authenticated version is available online at
Digital Object Identifier 10.1109/TSP.2012.2212436  
} \index{in
IEEE:\let\thefootnote\relax\footnote{Manuscript received December 7,
2012; revised May 30, 2012; accepted July 27, 2012. Copyright (c)
2012 IEEE. Personal use of this material is permitted. However,
permission to use this material for any other purposes must be
obtained from the IEEE by sending a request to
pubs-permissions\@ieee.org. The author is with  Department of
Mathematics and Statistics, Curtin University, GPO Box U1987, Perth,
Western Australia, 6845 (email N.Dokuchaev@curtin.edu.au). The
associate editor coordinating the review of this manuscript and
approving it for publication was Dr. Lawrence Carin. This work  was
supported by ARC grant of Australia DP120100928 to the author.}}
\\ {\rm Index terms} ---
Bandlimited, causal convolution, discrete time systems, harmonic
analysis,  prediction, Szeg\"o-Kolmogorov Theorem.
\\ MSC 2010 classification : 42A38, 
93E10, 
62M15,      
42B30,    
\end{abstract}
\def\break{}
\section{Introduction}
The paper studies pathwise predictability of discrete time processes
in a deterministic setting. It is well known that certain
restrictions on the frequency distribution can ensure additional
opportunities for prediction and interpolation of the processes. The
classical result for continuous time processes is the
Nyquist-Shannon-Kotelnikov interpolation theorem for the continuous
time band-limited processes.   These processes are presented in many
models including econometrics models (see examples in
\cite{pollock}). In theory, it is not possible to conclude that a
process is band-limited given some finite interval of observations.
In practice, this conclusion is being made based on historical data
for a certain process; this leads to models where processes are
assumed to be band-limited. Predictability based on sampling and the
Nyquist-Shannon-Kotelnikov theorem was discussed in \cite{W},
\cite{K}, \cite{M},\cite{V}, \cite{L} \cite{Ly}. These references
deal with the predictability of continuous-time, band-limited
stochastic processes, which include stationary processes. The
predictors obtained in this work were constructed for the setting
where the shape of the spectral representation is supposed to be
known. For discrete time processes, the predictability can also be
achieved given some properties of spectral representations. For
stationary discrete time Gaussian processes, Szeg\"o-Kolmogorov
Theorem ensures that the optimal prediction error is zero if the
spectral density $\phi$ is such that $\int_{-\pi}^\pi \log\phi\ew
d\o=-\infty$ (see, e.g., \cite{PM}, p.68). However, it was unknown
how to construct a predictor when  the shape of the spectral density
is unknown. This long standing problem was addressed in \cite{D12}:
a predictor for general type band-limited time series was suggested
in a deterministic setting. This predictor was a modification of the
predictor obtained in \cite{D08} for continuous time processes.
\par
Unfortunately, there are serious limitations to the practical  use
of the predictors constructed for the band-limited processes.  A
common argument dismissing the effectiveness of these  predictors is
that the predictors are not robust with respect to small noise
contamination. This leads to the conclusion  that the predictability
is an abnormality that disappears with the presence of arbitrarily
small noise or some incompleteness of historical data.
\par
This paper addresses these problems again. We consider discrete time
processes with some restrictions on the rate of energy decay on the
higher frequencies. We establish the predictability of these
processes, and  we suggest new linear predictors represented by
causal convolution sums over past times representing historical
observations  (Theorem \ref{ThM}).  The future values of the process
are not supposed to be calculated precisely but  are rather
approximated with an error in a prescribed interval that can be made
arbitrarily small uniformly over a wide class of underlying
processes.  Similarly to \cite{D12}, the predictors are given
explicitly in the frequency domain.  Whereas the results of
\cite{D12} were restricted to band-limited discrete-time processes,
we now extend the analysis to processes that are not strictly band
limited.  The predictors suggested here are different from the
predictors from \cite{D12}. The setting of the present paper is
similar to the one for continuous time processes from \cite{D10},
where predictors were suggested for processes that were not
band-limited but were assumed to have an exponential rate of energy
decay on higher frequencies.

These results  sheds some  new light on the predictability conundrum
for band-limited processes. More precisely, it leads to the
conclusion that the band-limited processes still allow  robust
predictability. The feasible predictability is an interval sense
only, i.e., it produces an interval that contains the future value
rather than the exact future value. The error can be made
arbitrarily small; however, this would require a large enough norm
value of the predictor's transfer function. We show that this
prediction is robust with respect to noise contamination in the
following sense. Given the size of the prediction error that is
associated with a process that is free of noise contamination, the
additional prediction error that is attributable to noise
contamination depends linearly on the product of the noise with a
norm associated with the transfer function entailed in forming the
prediction from preceding values of the process (see Section
\ref{secRob}).  If the predictor is targeting too small a size of
the error for processes without noise contamination, the norm of the
transfer function increases, and this robustness vanishes. \index{It
is also robust with respect to replacing the semi-infinite time
interval of observations by a large but finite one.}

The paper is organized in the following manner. In Section
\ref{secMain}, we formulate the definitions and the main result. In
Section \ref{secProof}, we prove the main theorem  concerning the
predictablility of processes with a certain rate of energy decay on
higher frequencies. In Section \ref{secRob}, we discuss the
robastness of the predictors with respect to noise contamination.
Finally, in Section \ref{secCon}, we summarize our results and offer
suggestions for further research.

\section{Definitions and main result}\label{secMain}
Let $D\defi\{z\in\C: |z|\le1\}$, $D^c=\C\backslash D$, and
$\T\defi\{z\in\C:\ |z|=1\}$.\index{ $\ZZ$ is the set of all
integers}
\par
We denote by $\ell_r$ the set of all sequences
$x=\{x(t)\}\subset\R$, $t=0,\pm 1,\pm 2$, such that
$\|x\|_{\ell_r}=\left(\sum_{t=-\infty}^{\infty}|x(t)|^r\right)^{1/r}<+\infty$
for $r\in[1,\infty)$ or  $\|x\|_{\ell_\infty}=\sup_t|x(t)|<+\infty$
for $r=+\infty$. \par Let $\ell_r^+$ be the set of all sequences
$x\in\ell_r$ such that $x(t)=0$ for $t=-1,-2,-3,...$.

\par
For  $x\in \ell_1$ or $x\in \ell_2$, we denote by $X=\Z x$ the
Z-transform  \baaa X(z)=\sum_{t=-\infty}^{\infty}x(t)z^{-t},\quad
z\in\C. \eaaa Respectively, the inverse $x=\Z^{-1}X$ is defined as
\baaa x(t)=\frac{1}{2\pi}\int_{-\pi}^\pi X\left(e^{i\o}\right)
e^{i\o t}d\o, \quad t=0,\pm 1,\pm 2,....\eaaa

 If $x\in \ell_2$, then $X|_\T$ is defined as an element of
$L_2(\T)$.
\par
Let  $H^r$ be the Hardy space of functions that are holomorphic on
$D^c$ including the point at infinity   (see, e.g.,  \cite{Du}
\index{Duren (1970)}). Note that Z-transform defines a bijection
between the sequences from $\ell_2^+$ and the restrictions (i.e.,
traces) of the functions from $H^2$ on $\T$.
\begin{definition}
Let $\w\K$ be the class of functions $\w k:\ZZ\to\R$
such that $\w k (t)=0$ for $t<0$ and such that $\w K(\cdot)=\Z\w
k\in H^\infty$.
\end{definition}
\begin{definition} Let $\Y\subset \ell_r$ be a class of
processes.
\begin{itemize}
\item[(i)] We say that this class is  $\ell_r$-predictable if
there exists a sequence $\{\w k_m(\cdot)\}_{m=1}^{+\infty}\subset
\w\K$ such that \baaa \|x(t+1)-\w x_m(t)\|_{\ell_r}\to 0\quad
\hbox{as}\quad m\to+\infty\quad\forall x\in\Y. \eaaa Here $  \w
x_m(t)\defi \sum^t_{s=-\infty}\w k_m(t-s)x(s).$
\item[(ii)]
 We say that the class $\Y$ is  uniformly $\ell_r$-predictable  if, for any $\e>0$, there exists $\w k(\cdot)\in \w\K$ such that \baaa \|x(t+1)- \w x(t)\|_{\ell_r}\le \e\quad
\forall x\in\Y. \label{predu}\eaaa Here $\w x(t)\defi
\sum^t_{s=-\infty}\w k(t-s)x(s).$
\end{itemize}
\end{definition}
\par
Let some $q>1$  be given. For $c>0$ and $\o\in[-\pi,\pi]$, set \baaa
h(\o,c)=\exp\frac{c}{[(\cos(\o)+1)^2+\sin^2(\o)]^{q/2}}.
\label{hdef}\eaaa
\par
Let $\X(c)$ be the class of all sequences $x\in\ell_2$ such that
\baa \esssup_{\o\in[-\pi,\pi]} |X\ew| h(\o,c)< +\infty,
\label{hfin}\eaa\ where $X=\Z x$. Let $\X=\cup_{c>0}\X(c)$.
\par
 Note that $h(\o,c)\to +\infty$ as $\o\to\pm \pi$ and that (\ref{hfin}) holds for degenerated processes, with $X\ew$
approaching zero with sufficient rate of decay as $\o\to \pm \pi$.
In particular, the class $\X$ includes all band-limited processes
$x$ such that $X\ew=0$ for $\o\notin[-\o_1,\o_1]$, for some
$\o_1\in(0,\pi)$, where $X=\Z x$.
\par

  \begin{theorem}\label{ThM}
Let either $r=2$ or $r=+\infty$. \begin{itemize}
\item[(i)]
The class $\X$ is  $\ell_r$-predictable.
\item[(ii)] Let $c_0>0$ be given, and let $\U(c_0)$   be a class of processes $x(\cdot)\in \X(c_0)$ such that
\baaa  \esssup_{\o\in[-\pi,\pi]} |X\ew h(\o,c_0)|\le 1\eaaa for
$X=\Z x$. \index{and that \baaa \int_{\o:\|\w|>\nu}|X\ew |^2d\o\to
0\quad\hbox{uniformly in}\quad x\in\U_0\quad \hbox{as}\quad \nu\to
0. \eaaa} Then this class $\U(c_0)$ is uniformly
$\ell_r$-predictable.
\item[(iii)]
A sequence of predicting kernels that ensures prediction required in
(i) and (ii) can be constructed as the following. Let $\mu>1$ be
given. For $\g>0$, set \baaa \a=\a(\g)=1-\g^{\frac{2\mu}{1-q}},\quad
V(z)\defi 1-\exp\left(-\frac{\g}{z+\a}\right),\brea\quad\w
K(z)\defi zV(z). \eaaa Then the required sequence of kernels that
ensures prediction required in (i) and (ii) is $\w k(\cdot)=\w
k(\cdot,\g)=\Z^{-1}\w K$, where $\g=\g_i\to +\infty$. For these
kernels, \baaa \|x(t+1)-\w x(t)\|_{\ell_r}\to 0\quad \hbox{as}\quad
\g\to +\infty\quad\forall x\in\X. \eaaa Moreover, for any $c_0>0$ and
$\e>0$, there exists $\g>0$ such that \baa \|x(t+1)- \w
x(t)\|_{\ell_r}\le \e\quad \forall x\in\U(c_0). \label{pred}\eaa
Here  $\w x(t)\defi \sum^t_{s=-\infty}\w k(t-s)x(s).$
\end{itemize}
\end{theorem}
Note that any particular  predictor described in Theorem \ref{ThM}
ensures predictability in an interval sense only, i.e., it produces
an interval $\left[\w x(t+1)-\e,\,\w x(t+1)+\e\right]$ that contains
 $x(t+1)$ rather than the exact value of $x(t+1)$. However, this $\e$
can done arbitrarily via selection of a large enough $\g$.
\par The family of
predicting kernels  $\w k$ introduced above  represents an extension
on the discrete time setting  of the construction introduced in
\cite{D10} for continuous time processes with exponential rate of
decay of energy on higher frequencies.
\section{Proofs}\label{secProof}
In our setting, $x(t+1)$ is the output of anticausal convolution
with the transfer function $K(z)\equiv z$, i.e., $x(t+1)=\Z^{-1}(K\Z
x)(t)$.
\par
 Let $\O(\a)=\arccos(-\a)$, let $D_+(\a)=(-\O(\a),\O(\a))$, and let
 a $D(\a)\defi[-\pi,\pi]\backslash D_+(\a)$. We have that  $\cos(\O(\a))+\a=0$, $\cos(\o)+\a>0$ for  $\o\in D_+(\a)$,
and $\cos(\o)+\a<0$ for  $\o\in D(\a)$.
\par
Note that $\a=\a(\g)\to 1$  as $\g\to+\infty$.
\begin{lemma}\label{lemmaV}
\begin{itemize}
\item[(i)] $V(z)\in H^{\infty}$ and $\w \Ko (z)\defi \Ko(z)V(z)\in
H^{\infty}$.
\item[(ii)] $V(e^{i\o})\to 1$ for all  $\o\in (-\pi,\pi)$ as  $\g\to +\infty$.
\item[(iii)] If $\o\in(-\O(\a),\O(\a))$  then $\Re \left(\frac{\g}{e^{i\o}+\a}\right) >0$ and $|V\ew -1|\le \rd{2}$.
\item[(iv)] For any $c>0$,
there exists $\g_0>0$ such that for any $\g\ge \g_0$ and for $V$
selected with $\a=\a(\g)$  we have $\int_{D(\a)}|V\ew-1|^\rho
h(\o,c)^{-\rho}d\o\le 2\arccos(\a)$ for any $\rho\ge 1$.
\end{itemize}
\end{lemma}
\par
{\it Proof of Lemma \ref{lemmaV}}. Clearly, $V\in H^{\infty}$, and
$zV(z)=K(z)V(z)\in H^{\infty}$, since the growth   of $z$  is being
compensated by multiplying with
$V(p)=-\sum_{k=1}^{+\infty}\frac{(-1)^k\g^k}{k!(z+\a)^k}$. Then statement (i)
follows.
\par
Further, for $\o\in(-\pi,\pi)$, we have that
 \baaa \g
\frac{1}{e^{i\o}+\a}=\g \frac{e^{-i\o}+\a}{|e^{i\o}+\a|^2}. \eaaa
Hence
 \baaa
\Re\left(\g \frac{1}{e^{i\o}+\a}\right)  =\g
\frac{\cos(\o)+\a}{|e^{i\o}+\a|^2}. \label{re} \eaaa If
$\g\to+\infty$ then $\a=\a(\g)\to 1$. This implies statements
(ii)-(iii).
\par
Let us prove statement (iv). We consider a large enough $\g$ such
that $\a=\a(\g)>3/4$. For these $\a$ and  $\o\in D(\a)$, we have
that $1/|e^{i\o}+\a|<2/|e^{i\o}+1|$. Hence
$\g|\Re(1/(e^{i\o}+\a))|<2\g/|e^{i\o}+1|$ for these $\a$ and $\o$.
Hence \baaa|V(i\o)-1|h(\o,c)^{-1}\le
\exp\left(\left|\Re\frac{\g}{e^{i\o}+\a}\right|-
\frac{c}{|e^{i\o}+1|^q}\right)\brea\le
\exp\left(\frac{2\g}{|e^{i\o}+1|}- \frac{c}{|e^{i\o}+1|^q}\right)
\eaaa for all $\o\in D(\a)$.
 Let $\rho(\a)=|e^{i\O(\a)}+1|^{-1}$. By the choice of $\a=\a(\g) $, it follows  that \baaa
 \rho(\a)=(2-2\a)^{-1/2}=\left(2\g^{\frac{2\mu}{1-q}}\right)^{-1/2}=2^{-1/2}\g^{\frac{\mu}{q-1}}.\eaaa
 Hence  \baaa
2\g
\rho(\a)^{1-q}=2\g\left(2^{-1/2}\g^{\frac{\mu}{q-1}}\right)^{1-q}=
2^{1-(1-q)/2}\g^{1-\mu}.\label{yaa}\eaaa It follows that \baa 2\g
\rho(\a)^{1-q}\to 0 \quad\hbox{as}\quad\g\to +\infty.\label{ya}\eaa
 By (\ref{ya}),
for any $c>0$, there exists $\g_0>0$ such that $2\g\rho(\a)^{1-q}\le
c$ for any $\g>\g_0$, and, therefore, $2\g \rho(\a)\le c\rho(\a)^q$.
Moreover, $2\g/|e^{i\o}+1|\le c/|e^{i\o}+1|^q$ for all $\o\in D(\a)$
since $1/|e^{i\o}+1|\ge \rho(\a)$ for these $\o$. Hence
\baaa|V(i\o)-1|h(\o,c)^{-1} \le 1,\quad |V(i\o)-1|^\rho
h(\o,c)^{-\rho} \le 1\eaaa for all $\o\in D(\a)$. In addition, the
measure of the set $D(\a)$ is
$\pi-\arccos(-\a)+\arccos(\a)=2\arccos(\a)$.
 This completes the proof of statement (iv) and Lemma \ref{lemmaV}. $\Box$
\vspace{0.5cm}
\par
{\it Proof of Theorem \ref{ThM}}. Let $\g\to +\infty$,  and let $V$,
$K$, $\w K$  be as defined above. Let $\ko=\Z^{-1}\Ko $ and $\w
\ko=\Z^{-1}\w \Ko$. For  $x(\cdot)\in \X$, let $X\defi \Z x$ and
\baaa &&\yo (t)\defi \sum_{s=t}^{\infty}\ko (t-s)x(s)=x(t+1),\quad \breakk \w
\yo (t)\defi \sum^t_{s=-\infty}\w \ko (t-s)x(s). \eaaa

We have that $\w k=\Z^{-1}\w K$ is real valued, since $k(\cdot)$ is
real valued and $K\left(\oo z\right)=\overline{K\left(z\right)}$,
$K\left(e^{-i\o}\right)=\overline{K\left(e^{i\o}\right)}$.

Let
 $\Yo \ew\defi (\Z \yo )\ew=K\ew X\ew$.
 By
the definitions, it follows that
 $\w Y\left(e^{i\o}\right)\defi \w K\left(e^{i\o}\right)X\left(e^{i\o}\right)=\Yo (\Z \w y )\ew$.
\par
Further, let  $\rho=2$ if $r=2$ and $\rho=1$ if $r=+\infty$.
\par
We have that $\|\w Y\ew-Y\ew\|_{L_\rho (-\pi,\pi)}^\rho=I_1+I_2,$
where \baaa &&I_1=\int_{D(\a)}|\w Y\ew-Y\ew|^\rho d\o,\qquad\breakk I_2=\int_{D_+(\a)}|\w Y\ew-Y\ew|^\rho d\o. \eaaa

By the assumptions, there exists $c>0$ such that $\|X\ew
h(\o,c)\|_{L_\infty(-\pi,\pi)}<+\infty$. Hence
\index{ \baaa I_1^{1/\rho}&=&\|\w
 Y\ew-Y\ew\|_{L_\rho (D(\a))}=
\|(\w K\ew-K\ew)X\|_{L_\rho (D(\a))}\nonumber\\
&\le& \|(V \ew -1)
 h(\o,c)^{-1}\|_{L_\rho (D(\a))} \|K\ew X\ew h(\o,c)\|_{L_\infty(-\pi,\pi)}\nonumber
\\&\le& (2\arccos(\a))^{1/\rho} \|X\ew h(\o,c)\|_{L_\infty(-\pi,\pi)} \rd{\|K\ew \|_{L_\infty(-\pi,\pi)}}.
\label{4s}\eaaa }
\rd {}
\baaa &&I_1^{1/\rho}\|\w
 Y\ew-Y\ew\|_{L_\rho (D(\a))}\breakk=
\|(\w K\ew-K\ew)X\|_{L_\rho (D(\a))}\nonumber\\
&&\le \|(V \ew -1)
 h(\o,c)^{-1}\|_{L_\rho (D(\a))}\breakk \times \|X\ew h(\o,c)\|_{L_\infty(-\pi,\pi)}\nonumber
\\&&\le (2\arccos(\a))^{1/\rho} \|X\ew h(\o,c)\|_{L_\infty(-\pi,\pi)}.
\label{4s}\eaaa

 The last inequality holds by  Lemma \ref{lemmaV}
(iv). It follows that $I_1\to 0$ as $\g\to +\infty$.

Let us estimate $I_2$. Lemma \ref{lemmaV} (iii) gives that
$|V\left(e^{i\o}\right)-1|\le \rd{2}$ for all $\o\in D_+(\a)$.
  We
 have that
\baaa I_2=\int_{D_+(\a)} |K\ew(1-V\ew) X\ew| ^\rho d\o\brea\le
\psi(\g)\|X\ew\|_{L_\infty(-\pi,\pi)}\\\le 2\psi(\g)\|X\ew
h(\o,c)\|_{L_\infty(-\pi,\pi)} ^\rho ,\eaaa where \baaa \psi(\g)
=\int_{D_+(\a)} |K\ew(1-V\ew)| ^\rho d\o\brea=\int_{-\pi}^\pi
\Ind_{D_+(\a)}(\o) |K\ew(1-V\ew)| ^\rho d\o. \eaaa Here $\Ind$
denotes the indicator function.
\par
  By Lemma \ref{lemmaV}(ii),
$\Ind_{D_+(\a)}(\o) |K\ew(1-V\ew)| ^\rho \to 0$ a.e. as $\g\to
+\infty$. By Lemma \ref{lemmaV}(iii), \baaa
\Ind_{D_+(\a)}(\o)|K\ew(1-V\ew) | ^\rho\brea\le \sup_{\o\in
D_+(\a)}\rd{|K\ew | ^\rho.}\eaaa   From Lebesgue Dominance Theorem,
it follows that $\psi(\g)\to 0$ as $\g\to+\infty$. It follows that
$I_1+I_2\to 0$ for any $c>0$, $x\in\X(c)$. By the definition of
$\rho$, we have that $1/\rho+1/r=1$. Hence $\|\w y-y\|_{\ell_r}\to
0$ as $\g\to +\infty$ for any $x\in\X$.
 This
completes the proof of statement (i).

Let us prove statement (ii). We have that \baaa \|\w
Y\ew-Y\ew\|_{L_\rho (-\pi,\pi)} ^\rho = I_1+I_2\brea\le (2
\arccos(\a)+ 2\psi(\g))\|X\ew h(\o,c_0)\| ^\rho
_{L_\infty(-\pi,\pi)}\\\le 2 \arccos(\a)+ 2\psi(\g)\eaaa for any
$x\in\U(c_0)$. For any $\e>0$, one can select $\g$ such that
$2\psi(\g)\le \e ^\rho /2$ and that $2 \arccos(\a)\le \e^\rho/2$.
This choice ensures that $\|\w y-y\|_{\ell_r}\le \e$.  This
completes the proof of statement (ii). It follows that the
predicting kernels $\w k(\cdot)=\Z^{-1}\w K$ are such as required.
This completes the proof of Theorem \ref{ThM}. $\Box$
\par
It can be noted that the choice of predicting kernels is not unique.
In particular, the kernels preserve the  properties described in
Theorem \ref{ThM} for any selection of $\a=\a(\g)$ such that
(\ref{ya}) holds. For instance, $\a(\g)$ can be selected as \baaa
\a=\a(\g)=1-(\log \g)^{-1}\g^{\frac{2}{1-q}}. \eaaa In addition, it
follows from the proofs that the uniform predictability from
statement (ii) of Theorem \ref{ThM} can be ensured with
$$
\a=\a(\g)=1-\frac{1}{2}\left(2\g/c_0\right)^{\frac{2}{1-q}}.
$$ In this case,
$\rho(\a)=(2\g/c_0)^{1/(q-1)}$.
 This corresponds to the case where
$\mu=1$ in Theorem \ref{ThM}(iii).
\def\NN{{\scriptscriptstyle N}}
\section{On the prediction error generated by noise
contamination}\label{secRob} Let us estimate the prediction error
for the case when the predictor designed for processes from $\X$ is
applied to a process with a small high-frequency noise
contamination. Let us consider a process $x(\cdot)\in\ell_{\infty}$
such that $x=x_0+x_{\NN}$, where $x_0\in\X$,
$x_{\NN}\in\ell_{\infty}$. The process $x_{\NN}$ represents the
noise. Let $X=\Z x$, $X_0=\Z x_0$, and $X_{\NN}=\Z x_{\NN}$. We
assume that $X_0\ew\in L_1(-\pi,\pi)$ and
 $\|X_{\NN}\ew\|_{L_1(-\pi,\pi)}=\nu$. The parameter $\nu\ge 0$ represents the
intensity of the noise.
\par
Assume that the predictor  is constructed as in Theorem \ref{ThM}
under the hypothesis that $\nu=0$ (i.e, that $x_{\NN}=0$ and
$x\in\X$). For an arbitrarily small $\e>0$, there exists $\g$ such
that, if the hypothesis that $\nu=0$ is correct, then \baaa
\int_{-\pi}^{\pi}|(\w K\ew-K\ew
)X\ew|d\o\brea=\int_{-\pi}^{\pi}|(\w K\ew-K\ew)X_0\ew|d\o\le
2\pi\e\eaaa and\baaa \|\w
y-y\|_{\ell_{\infty}}\le\e.\label{eps}\eaaa
\par
Let us estimate the prediction error for the case where $\nu>0$. We
have that \baaa \|\w y-y\|_{\ell_{\infty}}\le J_0+ J_{\NN},\eaaa
where\baaa J_0=\frac{1}{2\pi}\|(\w
K\ew-K\ew)X_0\ew|\|_{L_1(-\pi,\pi)},\brea\quad
J_{\NN}=\frac{1}{2\pi}\|(\w K\ew-K\ew)X_{\NN}\ew|\|_{L_1(-\pi,\pi)}.
\eaaa The value $J_{\NN}$ represents the additional error caused by
the presence of unexpected high-frequency noise (when $\nu>0$). It
follows that \baa \|\w y-y\|_{\ell_{\infty}}\le
\e+\nu(\kappa+1),\label{yn}\eaa  where
$\kappa=\sup_{\o\in[-\pi,\pi]}|\w K\ew|$.
\par
Therefore, it can be concluded that the prediction  is robust with
respect to noise contamination for any given $\e$. On the other
hand, if $\e\to 0$ then $\g\to +\infty$ and $\kappa\to +\infty$. In
this case, error (\ref{yn}) is increasing for any given $\nu>0$.
This happens when the predictor is targeting too small a size of the
error for the processes from $\X$, i.e., under the assumption that
$\nu=0$.
\par
The equations describing the dependence of $\e$ and $\k$ on $\g$
could be derived similarly to estimates in \cite{D12}, Section 6,
where it was done for different predicting kernels and for
band-limited processes. We leave it for future research.
\section{Concluding remarks}\label{secCon}
\begin{itemize}
\item[(i)]
Technically, the predictors obtained above require the past values
of $x(s)$ for all $s\in(-\infty,t]$. However, $\sum_{s=-\infty}^t\w
k(t-s)x(s) $ can be approximated by $\sum_{s=-M}^t\w k(t-s)x(s) $
for a large enough $M>0$.  Therefore, the predictors are robust with
respect to replacing the semi-infinite time interval of observations
by a large finite one.
\item[(ii)]
 For processes from $\X$, the selection of a large enough $\g$
can ensure that the prediction  error is arbitrarily small. However,
 the corresponding  predictor transfer function $\w K$ is
large in norm for large $\g$. This leads to a large error caused by
any noise contamination presented for processes $x(\cdot)\notin\X$,
i.e., with the energy on higher frequencies that is not decaying
fast enough near the point $z=e^{i\pi}$. Nevertheless, the suggested
predictors are robust with respect to the contamination  noise for
any fixed $\g$. The error generated by the noise is limited if
$\kappa$ in (\ref{yn}) is limited, i.e., if $\g$ is limited and $\e$
in (\ref{yn}) is not too small.
\item[(iii)] The presence of robustness mentioned above leads to the
conclusion that certain interval type predictability for
band-limited processes is not an abnormality that disappears with
the presence of arbitrarily small noise or some incompleteness of
historical data. This predictability holds even for processes with a
certain rate of energy decay on higher frequencies that are not
exactly band-limited. This must be taken into account for all models
where band-limited processes are assumed. In particular, this
implies that band-limited processes should be used with caution in
the models where predictability is difficult to justify such as
models for financial time series.
\item[(iv)] The results of this paper can be applied to discrete time
stationary random Gaussian processes with the spectral densities
$\phi$ such that\index{$\int_{-\pi}^\pi \log\phi\ew d\o=-\infty$}
$\int_{-\pi}^\pi \phi\ew h(\o,c) d\o<+\infty$, i.e., when the
spectral density is decaying fast enough on the higher frequencies.
By Szeg\"o-Kolmogorov Theorem, it was known in principle that the
minimal (optimal) predicting error is zero in this case. However, it
was unknown how to construct the corresponding predictors  for
general classes of $\phi$.
\item[(v)] The restrictions on the spectral representations imposed
on the underlaying processes are quite tight. They are not satisfied
for the samples generated from autoregressive stochastic models
including ARMA, FARIMA, FEXP and other common long-memory processes
(see, e.g., \cite{R}). It is yet unclear in which applications the
processes with the required rate of energy decay on higher
frequencies could be found. Therefore, it is unclear where to find
real data sets to test the suggested predictors.  We may suggest
applying the corresponding predictors for band-limited processes
such as described in \cite{pollock}. We leave it for future
research.
\end{itemize}
\subsection*{Acknowledgment} This work  was
supported by ARC grant of Australia DP120100928 to the author. The
author thanks Prof. Augusto Ferrante for useful discussion and
advice on the time series analysis. \index{IEEE: Thanks are also due
to the anonymous reviewers for their comments and constructive
criticism.}

\end{document}